\documentclass[reprint,a4paper,11pt,reqno]{amsart}

\usepackage{amsmath}    
\usepackage{amsfonts}
\usepackage{amscd}
\usepackage[latin2]{inputenc}
\usepackage{t1enc}
\usepackage[mathscr]{eucal}
\usepackage{indentfirst}
\usepackage{graphicx}
\usepackage{graphics}
\usepackage{pict2e}
\usepackage{fancyhdr}
\numberwithin{equation}{section}
\usepackage[margin=2.9cm]{geometry}
\usepackage{epstopdf} 
\usepackage{hyperref}
\usepackage{inputenc}
\usepackage{mathtools}

\theoremstyle{plain}
\newtheorem{Th}{Theorem}
\newtheorem{Lemma}[Th]{Lemma}

\theoremstyle{definition}

\newtheorem{?}[Th]{Problem}

\hypersetup{
    colorlinks=true,
    linkcolor=blue,
    citecolor=blue,
    filecolor=magenta,      
    urlcolor=blue,
}

\begin{document}

\pagestyle{fancy}
\lhead{}
\chead{}
\rhead{\thepage}
\cfoot{}
\lfoot{}
\rfoot{}
\renewcommand{\headrule}{}

\title{A sequence of elementary integrals related to integrals studied by Glaisher that contain trigonometric and hyperbolic functions} 

\author{Martin Nicholson} 

\begin{abstract}  
We generalize several integrals studied by Glaisher. These ideas are then applied to obtain an analog of an integral due to Ismail and Valent. 
\end{abstract}

\maketitle

\section{Introduction}

\thispagestyle{empty}

The following integral
\begin{equation}\label{glaisher1}
    \int\limits_0^\infty \frac{\sin x\sinh (x/a)}{\cos (2 x)+\cosh \left(2x/a\right)}\frac{dx}{x}=\frac{\tan^{-1} a}{2}
\end{equation}
can be deduced as a particular case of entry $4.123.6$ from \cite{GR}. The case $a=1$ of this integral can be found in an old paper by Glaisher \cite{glaisher}. More recently, symmetric cases were also stated in \cite{glasser} and \cite{pm}. We are going to generalize the above integral as
\begin{Th}\label{th1} Let $n$ be an odd integer, then
\begin{align}\label{form1}
\int\limits_0^{1}\frac{\sin \bigl(n \sin^{-1}t\bigr)\sinh \bigl(n \sinh^{-1}(t/a)\bigr)}{\cos \bigl( 2 n \sin^{-1}t\bigr)+\cosh \bigl(2 n \sinh^{-1}(t/a)\bigr)}\frac{dt}{t \sqrt{1-t^2} \sqrt{1+{t^2}/{a^2}}}=\frac{\tan^{-1} a}{2}.
\end{align}
\end{Th}
\noindent When $n$ is large, then the main contribution to the integral \ref{form1} comes from a small neighbourhood around $t=0$ and the integral reduces to \ref{glaisher1}.

Another integral by Glaisher reads (equation $24$ in \cite{glaisher})
\begin{equation}\label{glaisher2}
    \int\limits_0^\infty \frac{\cos x\cosh x}{\cos (2 x)+\cosh \left(2x\right)}\,x\,dx=0.
\end{equation}
It would be generalized as 

\begin{Th}\label{th2} Let $n$ be an even integer, then
\begin{equation}\label{}
\int\limits_0^{1}\frac{\cos (n \sin^{-1}t)\cosh \left(n \sinh^{-1}t\right)}{\cos (2 n \sin^{-1}t)+\cosh \left(2 n \sinh^{-1}t\right)}\frac{tdt}{\sqrt{1-t^4} }=0.
\end{equation}
\end{Th}
\noindent Unfortunately there doesn't seem to be any nice parametric extensions similar to that in Theorem \ref{th1}.

A particularly interesting integral is 
\[
\int\limits_{-\infty}^{\infty}\frac{dt}{\cos (K\sqrt{t})+\cosh \left(K'\sqrt{t}\right)}=1,
\]
studied by Ismail and Valent \cite{ismail}. Here $K=K(k)$ and $K'=K(\sqrt{1-k^2})$ are elliptic integrals of the first kind. 
Berndt \cite{berndt} gives a generalization of this formula and as an intermediate result proves that (see Corollary $3.3$)
\begin{equation}\label{berndt1}
    \int\limits_0^\infty\frac{x^{4s+1}dx}{\cos x+\cosh x}=(-1)^s\frac{\pi^{4s+2}}{2^{2s+1}}\sum_{j=0}^\infty (-1)^{j}\frac{(2j+1)^{4s+1}}{\cosh\frac{\pi  (2 j+1)}{2 n}}
\end{equation}
for non-negative integer $s$. The next theorem gives an elementary analog of \ref{berndt1}.
\begin{Th}\label{th3} Let $s$ and $n$ be a integers such that $0\le s<\left\lfloor{\frac{n}{2}}\right\rfloor$. Then
\begin{align}
\nonumber\int\limits_0^{1}\frac{t^{2s}}{\cos (2 n \sin^{-1}\!\sqrt{t})+\cosh \left(2 n \sinh^{-1}\!\sqrt{t}\right)}\frac{dt}{\sqrt{1-t^2} }\phantom{...............................}\\=\frac{\pi  (-1)^s}{2^{2 s+1} n}\sum _{j=1}^{n/2} \frac{(-1)^{j-1} \tan \frac{\pi  (2 j-1)}{2 n}}{\cosh \left(n \sinh ^{-1}\tan \frac{\pi  (2 j-1)}{2 n}\right)}\left(\frac{\sin ^2\frac{\pi  (2 j-1)}{2 n}}{\cos \frac{\pi  (2 j-1)}{2 n}}\right)^{2 s}.
\end{align}
\end{Th}
Kuznetsov \cite{kuznetsov} proved that
\begin{equation}
\frac{1}{2}\int_{\mathbb{R}}  \frac{\sin(\sqrt{x}u)}{{\sqrt{x}}} \cdot\frac{dx}
{\cos(\sqrt{x}K)+\cosh(\sqrt{x}K')}=
 \frac{{\textnormal{sn}}(u,k){\textnormal{dn}}(u,k)}{{\textnormal{cn}}(u,k)}.
\end{equation}
By differentiating with respect to $u$ one can deduce
\begin{equation}
\frac{1}{2}\int_{\mathbb{R}} \frac{\cos(\sqrt{x}u)}
{\cos(\sqrt{x}K)+\cosh(\sqrt{x}K')}\, dx=k^2{\textnormal{cn}^2}(u,k)+
 \frac{1-k^2}{{\textnormal{cn}^2}(u,k)}.
\end{equation}
Right hand side of this formula has a Fourier series expansion which can be obtained using fundamental relations and the expansion of ${\textnormal{sn}^2}(u,k)$ given in \cite{ww}:
\begin{equation}\label{fourier}
    \frac{\pi^2}{4K^2\cos^2\frac{\pi u}{2K}}+\frac{\pi^2}{K^2}\sum _{j=1}^\infty j(-1)^{j-1}\left(\left\{\tanh\frac{\pi  jK'}{2 K}\right\} ^{(-1)^j}-1\right)\cdot \cos\frac{\pi ju}{K}.
\end{equation}
Stated in this way, symmetric case $K=K'$ of Kuznetsov's formula admits a finite analog of the form
\begin{Th}\label{th4} Let $n$ and $u$ be integers such that $|u|<n$. Then
\begin{align*}
\nonumber\int\limits_{-1}^1&\frac{\cos \left(2 u \sin ^{-1}\sqrt{t}\right)}{\cos \left(2 n \sin ^{-1}\sqrt{t}\right)+\cosh \left(2 n \sinh ^{-1}\sqrt{t}\right)}\frac{dt}{\sqrt{1-t^2} }\\&\phantom{........}=\frac{\pi}{2n}\sum_{j=1}^{2n}\frac{ (-1)^{j-1}\sin\frac{\pi  j}{2 n} }{\sqrt{1+\sin ^2\frac{\pi j}{2 n}}}\left\{\tanh\left(n \sinh ^{-1}\sin\frac{\pi  j}{2 n}\right)\right\} ^{(-1)^j}\!\!\!\!\cdot\cos\frac{\pi j u}{n}
\end{align*}
\end{Th}
\noindent The finite analog of the term with $\sec^2\frac{\pi u}{2K}$ in \ref{fourier} is accounted for by the sum valid for integer $u$
\[
\sum _{j=1}^{2 n} (-1)^{j-1} \sin \frac{\pi  j}{2 n}\cos\frac{\pi  ju}{n}=\frac{\sin \frac{\pi }{2 n}}{\cos \frac{\pi }{2 n}+\cos\frac{\pi u}{n}}.
\]

Proofs of these theorems are given in the subsequent sections \ref{proof1}, \ref{proof2}, \ref{proof3}, and  \ref{proof4}. In section \ref{disc}, some discussions of the theorems are give. In particular, it will be explained that the form of the integral in Theorem \ref{th3} is not arbitrary. Its form has been chosen to reflect a certain kind of symmetry satisfied also by integrals in Theorems \ref{th1} and \ref{th2}. Some open questions will be discussed in section \ref{open}.

\section{Proof of Theorem \ref{th1}}\label{proof1}
We break the proof into a series of lemmas.
\begin{Lemma}\label{lemma1} Let $n$ be an odd integer. Then we have the partial fractions expansion
\begin{align}\label{ae1}
\nonumber\frac{\sin \bigl(n \sin^{-1}t\bigr)\sinh \bigl(n \sinh^{-1}(t/a)\bigr)}{\cos \bigl( 2 n \sin^{-1}t\bigr)+\cosh \bigl(2 n \sinh^{-1}(t/a)\bigr)}\frac{2n}{t^2}\phantom{.............................................}\\=\sum _{j=1}^n\frac{i(-1)^{j-1} }{\sin\frac{\pi  (2 j-1)}{2 n}}\cdot \frac{\left(a\cos\frac{\pi  (2 j-1)}{2 n}+i\right) \left(a+i \cos\frac{\pi  (2 j-1)}{2 n}\right)}{t^2 \left(a^2-1+2 ia \cos\frac{\pi  (2 j-1)}{2 n}\right)-a^2 \sin ^2\frac{\pi  (2 j-1)}{2 n}}.
\end{align}
\end{Lemma}

\noindent{\it{Proof.}} When $n$ is an odd integer, the expressions 
\[
2n\sin \bigl(n \sin^{-1}t\bigr)\sinh \bigl(n \sinh^{-1}(t/a)\bigr)/t^2,\quad \cos \bigl( 2 n \sin^{-1}t\bigr)+\cosh \bigl(2 n \sinh^{-1}(t/a)\bigr)
\]
are polynomials in $t^2$ of degrees $n-1$ and $n$, respectively:
\[
\frac{\sin \bigl(n \sin^{-1}t\bigr)\sinh \bigl(n \sinh^{-1}(t/a)\bigr)}{\cos \bigl( 2 n \sin^{-1}t\bigr)+\cosh \bigl(2 n \sinh^{-1}(t/a)\bigr)}\frac{2n}{t^2}=\frac{P_{n-1}(t^2)}{Q_n(t^2)}.
\]

Let us find the $n$ roots of the denominator polynomial $Q_n(x)$. $Q_n(x)$ can be written as 
\[
Q_n(x)=\cos \bigl( n \sin^{-1}\sqrt{x}+in \sinh^{-1}(\sqrt{x}/a)\bigr)\cos \bigl( n \sin^{-1}\sqrt{x}-in \sinh^{-1}(\sqrt{x}/a)\bigr),
\]
and thus its roots can be found from the equations
\[
\sin^{-1}\sqrt{x}\pm i\sinh^{-1}(\sqrt{x}/a)=\frac{\pi(2j-1)}{2n}, \quad j=1,2,..., n,
\]
or equivalently from the equations
\[
\sqrt{x}\sqrt{1+\frac{x}{a^2}}\pm \frac{i\sqrt{x}}{a}\sqrt{1-x}=\sin\frac{\pi(2j-1)}{2n}, \quad j=1,2,..., n,
\]
One can get rid of the radicals to come to a quadratic equation with respect to $x$:
\[
x^2\left((1-a^2)^2+4a^2\cos^2\frac{\pi(2j-1)}{2n}\right)+2xa^2(1-a^2)\sin^2\frac{\pi(2j-1)}{2n}+\sin^4\frac{\pi(2j-1)}{2n}=0, \quad j=1,2,..., n.
\]
One can easily deduce from this that the $n$ roots of the denominator polynomial are
\begin{equation}\label{roots}
    x_j=\left(a^2-1+2 ia \cos\frac{\pi  (2 j-1)}{2 n}\right)^{-1}a^2 \sin ^2\frac{\pi  (2 j-1)}{2 n}, \quad j=1,2,..., n.
\end{equation}

Now we can find the partial fractions expansion
\begin{equation}\label{partfrac}
    \frac{P_{n-1}(t^2)}{Q_n(t^2)}=\sum_{j=1}^n\frac{P_{n-1}(x_j)}{Q_n^{\,'}(x_j)}\frac{1}{t^2-x_j}.
\end{equation}
A simple calculation shows that
\[
\frac{Q_n^{\,'}(x_j)}{P_{n-1}(x_j)}=\sqrt{\frac{x_j}{a^2+x_j}}\frac{\cosh\bigl(n \sinh^{-1}(\sqrt{x_j}/a)\bigr)}{\sin\bigl(n \sin^{-1}\sqrt{x_j}\bigr)}-\sqrt{\frac{x_j}{1-x_j}}\frac{\cos\bigl(n \sin^{-1}\sqrt{x_j}\bigr)}{\sinh\bigl(n \sinh^{-1}(\sqrt{x_j}/a)\bigr)}.
\]
The equation $\cos \bigl( 2 n \sin^{-1}\sqrt{x_j}\bigr)+\cosh \bigl(2 n \sinh^{-1}(\sqrt{x_j}/a)\bigr)=0$ implies
\[
\cosh\bigl(n \sinh^{-1}(\sqrt{x_j}/a)\bigr)=\mu_j\sin\bigl(n \sin^{-1}\sqrt{x_j}\bigr),\qquad \cos\bigl(n \sin^{-1}\sqrt{x_j}\bigr)=i\nu_j\sinh\bigl(n \sinh^{-1}(\sqrt{x_j}/a)\bigr),
\]
where $\mu_j=\pm,~\nu_j=\pm$. To determine the signs $\mu_j,\nu_j$, one can consider the limiting case $a>>1$. We have
\[
\sqrt{x_j}=\sin\frac{\pi  (2 j-1)}{2 n}-\frac{i}{a}\sin\frac{\pi  (2 j-1)}{2 n}\cos\frac{\pi  (2 j-1)}{2 n}+O(a^{-2}).
\]
This means
\[
\sin^{-1}\sqrt{x_j}=\frac{\pi  (2 j-1)}{2 n}-\frac{i}{a}\sin\frac{\pi  (2 j-1)}{2 n}+O(a^{-2}).
\]
From this it follows that $\mu_j=\nu_j=(-1)^{j-1}$ and thus
\begin{align*}
    \frac{Q_n^{\,'}(x_j)}{P_{n-1}(x_j)}&=(-1)^{j-1}\left(\sqrt{\frac{x_j}{a^2+x_j}}-i\sqrt{\frac{x_j}{1-x_j}}\right)\\
    &=i(-1)^j\frac{\sin \frac{\pi  (2 j-1)}{2 n}\left(a^2-1+2 ia \cos\frac{\pi  (2 j-1)}{2 n}\right)}{\left(a\cos\frac{\pi  (2 j-1)}{2 n}+i\right) \left(a+i \cos\frac{\pi  (2 j-1)}{2 n}\right)}.
\end{align*}
Substituting this into \ref{partfrac} we get the desired result. \qed

\begin{Lemma}\label{lemma2}
\begin{align}
\nonumber\int\limits_0^1 \frac{1}{t^2 \left(a^2-1+2 ia \cos\frac{\pi  (2 j-1)}{2 n}\right)-a^2 \sin ^2\frac{\pi  (2 j-1)}{2 n}}\frac{t\,dt}{\sqrt{1-t^2} \sqrt{1+{t^2}/{a^2}}}\\=\frac{\tan^{-1}a+i\tanh^{-1}\cos\frac{\pi  (2 j-1)}{2 n}}{i\left(a\cos\frac{\pi  (2 j-1)}{2 n}+i\right) \left(a+i \cos\frac{\pi  (2 j-1)}{2 n}\right)}.
\end{align}
\end{Lemma}

\noindent{\it{Proof.}} Composition of two substitutions $t^2=1-(1+1/a^2)\sin^2\phi$, ($0<\phi<\tan^{-1}a$) and $\tan\phi=s$, ($0<s<a$) reduces this integral to an integral of a rational function.\qed

\begin{Lemma}\label{lemma3} For $n$ odd, one has
\[
\sum _{j=1}^n \frac{(-1)^{j-1}}{\sin \frac{\pi  (2 j-1)}{2 n}}=n.
\]
\end{Lemma}

\noindent{\it{Proof.}} Put $t=1$, $a=i$ in Lemma \ref{lemma1}.\qed

From the three lemmas above it follows immediately that
\begin{align*}\label{form1}
\int\limits_0^{1}\frac{\sin \bigl(n \sin^{-1}t\bigr)\sinh \bigl(n \sinh^{-1}(t/a)\bigr)}{\cos \bigl( 2 n \sin^{-1}t\bigr)+\cosh \bigl(2 n \sinh^{-1}(t/a)\bigr)}\frac{dt}{t \sqrt{1-t^2} \sqrt{1+{t^2}/{a^2}}}\phantom{......}\\=\frac{\tan^{-1} a}{2}+\frac{i}{2n}\sum_{j=1}^n\frac{(-1)^{j-1}}{\sin\frac{\pi  (2 j-1)}{2 n}}\tanh^{-1}\cos\frac{\pi  (2 j-1)}{2 n}.
\end{align*}
To finish the proof, note that the sum in this formula is $0$ because (since $n$ is odd) $j$-th and $(n+1-j)$-th terms cancel each other out.

\section{Proof of Theorem \ref{th2}}\label{proof2}
\begin{Lemma}\label{lemma4} Let $n$ be an even integer. Then
\begin{align*}
&\frac{\cos (n \sin^{-1}t)\cosh \left(n \sinh^{-1}(t/a)\right)}{\cos (2 n \sin^{-1}t)+\cosh \left(2 n \sinh^{-1}(t/a)\right)}=\frac{(-1)^{n/2}}{2}\frac{a^n}{1+a^{2n}}\\&\phantom{.................}+\sum _{j=1}^{n}\frac{(-1)^{j} a^2\sin\frac{\pi  (2 j-1)}{2 n}}{2n\left(a^2-1+2 ia \cos\frac{\pi  (2 j-1)}{2 n}\right)}\cdot \frac{\left(a\cos\frac{\pi  (2 j-1)}{2 n}+i\right) \left(a+i \cos\frac{\pi  (2 j-1)}{2 n}\right)}{t^2 \left(a^2-1+2 ia \cos\frac{\pi  (2 j-1)}{2 n}\right)-a^2 \sin ^2\frac{\pi  (2 j-1)}{2 n}}.
\end{align*}
\end{Lemma}
\noindent{\it{Proof.}} When $n$ is even, the functions $\cos (n \sin^{-1}t)$ and $\cosh \left(n \sinh^{-1}(t/a)\right)$ are a polynomials in $t^2$ of degree $n/2$. This means we can write
\[
\frac{\cos (n \sin^{-1}t)\cosh \left(n \sinh^{-1}(t/a)\right)}{\cos (2 n \sin^{-1}t)+\cosh \left(2 n \sinh^{-1}(t/a)\right)}=C+\frac{R_{n-1}(t^2)}{Q_n(t^2)},
\]
where $R_{n-1}$ is a polynomial of order $n-1$ and $Q_n$ was defined in the proof of the Lemma \ref{lemma1}. $Q_n(x)$ has $n$ roots given by \ref{roots}.

To find the constant $C$ consider the limit $t\to+\infty$ assuming that $a>0$. In this case
\[
\sin^{-1}t=\frac{\pi}{2}-i\ln(2t)+O(t^{-1}), \qquad \sinh^{-1}(t/a)=\ln(2t/a)+O(t^{-1}),
\]
and we get
\[
C=\frac{(-1)^{n/2}}{2}\frac{a^n}{1+a^{2n}}.
\]
Since the order of the polynomial $R_{n-1}$ is smaller than the order of the polynomial $Q_n$ we can write the partial fractions expansion
\[
\frac{R_{n-1}(t^2)}{Q_n(t^2)}=\sum_{j=1}^n\frac{R_{n-1}(x_j)}{Q_n^{\,'}(x_j)}\frac{1}{t^2-x_j}.
\]
A calculation similar to that in Lemma \ref{lemma1} shows that
\begin{align*}
\frac{Q_n^{\,'}(x_j)}{R_{n-1}(x_j)}&=\frac{2n}{x_j}\left(\sqrt{\frac{x_j}{a^2+x_j}}\frac{\sinh\bigl(n \sinh^{-1}(\sqrt{x_j}/a)\bigr)}{\cos\bigl(n \sin^{-1}\sqrt{x_j}\bigr)}-\sqrt{\frac{x_j}{1-x_j}}\frac{\sin\bigl(n \sin^{-1}\sqrt{x_j}\bigr)}{\cosh\bigl(n \sinh^{-1}(\sqrt{x_j}/a)\bigr)}\right)\\
&=(-1)^{j-1}\frac{2n}{x_j}\left(\sqrt{\frac{x_j}{a^2+x_j}}-i\sqrt{\frac{x_j}{1-x_j}}\right)\\
&=\frac{2n(-1)^{j}\left(a^2-1+2 ia \cos\frac{\pi  (2 j-1)}{2 n}\right)^2}{a^2\sin\frac{\pi  (2 j-1)}{2 n}\left(a\cos\frac{\pi  (2 j-1)}{2 n}+i\right) \left(a+i \cos\frac{\pi  (2 j-1)}{2 n}\right)}.
\end{align*}
This completes the proof of the lemma. \qed

Using Lemmas \ref{lemma2} and \ref{lemma4} we find
\begin{align*}
    &\int\limits_0^{1}\frac{\cos \bigl(n \sin^{-1}t\bigr)\cosh \bigl(n \cosh^{-1}(t/a)\bigr)}{\cos \bigl( 2 n \sin^{-1}t\bigr)+\cosh \bigl(2 n \sinh^{-1}(t/a)\bigr)}\frac{t\,dt}{\sqrt{1-t^2} \sqrt{1+{t^2}/{a^2}}}
    \\&=\frac{(-1)^{n/2}}{2}\frac{a^{n+1}}{1+a^{2n}}\tan^{-1}(1/a)+a^2\sum_{j=1}^n(-1)^{j}\,\frac{\tanh^{-1}\cos\frac{\pi  (2 j-1)}{2 n}-i\tan^{-1}a}{2n\left(a^2-1+2 ia \cos\frac{\pi  (2 j-1)}{2 n}\right)}\,\sin\frac{\pi  (2 j-1)}{2 n},
\end{align*}
and in particular when $a=1$
\begin{align*}
\int\limits_0^{1}\frac{\cos (n \sin^{-1}t)\cosh \left(n \sinh^{-1}t\right)}{\cos (2 n \sin^{-1}t)+\cosh \left(2 n \sinh^{-1}t\right)}\frac{tdt}{\sqrt{1-t^4} }&=\frac{\pi}{16}(-1)^{n/2}-\sum_{j=1}^n(-1)^{j}\,\frac{\pi+4i\tanh^{-1}\cos\frac{\pi  (2 j-1)}{2 n}}{16n\cot\frac{\pi  (2 j-1)}{2 n}}\\
&=\frac{\pi}{16n}\Bigg((-1)^{n/2}n-\sum_{j=1}^n {(-1)^{j}}\tan\frac{\pi  (2 j-1)}{2 n}\Bigg).
\end{align*}
To calculate the sum in this expression we use Lemma \ref{lemma4} with $t=1$ and $a\to\infty$ to get
\[
\sum_{j=1}^n{(-1)^{j}}\tan\frac{\pi  (2 j-1)}{2 n}=(-1)^{n/2}n.
\]
This completes the proof of the theorem.

\section{Proof of Theorem \ref{th3}}\label{proof3}

Here we restrict the consideration to the symmetric case $a=1$.
\begin{Lemma}\label{lemma5} The following partial fractions expansion holds for integers $s$ and $n$ such that $0\le s<\left\lfloor{\frac{n}{2}}\right\rfloor$
\begin{align*}
    &\frac{t^{2 s}}{\cos \bigl(2 n \sin ^{-1}\sqrt{t}\bigr)+\cosh \bigl(2 n \sinh ^{-1}\sqrt{t}\bigr)}\\&\phantom{.............}=\frac{(-1)^s}{2^{2 s} n}\sum_{j=1}^{n/2}\frac{1}{4t^2+\frac{\sin ^4\frac{\pi  (2 j-1)}{2 n}}{\cos ^2\frac{\pi  (2 j-1)}{2 n}}}\frac{(-1)^{j-1} \tan\frac{\pi  (2 j-1)}{2 n}}{\cosh \left(n \sinh ^{-1}\tan \frac{\pi  (2 j-1)}{2 n}\right)}\frac{1+\cos ^2\frac{\pi  (2 j-1)}{2 n}}{\cot ^2\frac{\pi  (2 j-1)}{2 n}}\left(\frac{\sin ^2\frac{\pi  (2 j-1)}{2 n}}{\cos\frac{\pi  (2 j-1)}{2 n}}\right)^{2s}.
\end{align*}
\end{Lemma}

\noindent{\it{Proof.}} From consideration of the limit $t\to +\infty$ once can see (similarly to that in Lemma \ref{lemma4}) that the leading coefficient of the polynomial $Q_n(t)=\cos \bigl(2 n \sin ^{-1}\sqrt{t}\bigr)+\cosh \bigl(2 n \sinh ^{-1}\sqrt{t}\bigr)$ is $2^{2n-1}(1+(-1)^n)$ and thus that $Q_n(t)$ is an even  polynomial of degree $2\left\lfloor{\frac{n}{2}}\right\rfloor$. Its roots are (see \ref{roots})
\[
x_j=-\frac{i\sin ^2\frac{\pi  (2 j-1)}{2 n}}{2\cos\frac{\pi  (2 j-1)}{2 n}},\quad y_j=\frac{i\sin ^2\frac{\pi  (2 j-1)}{2 n}}{2\cos \frac{\pi  (2 j-1)}{2 n}},\quad j=1,2,...,\left\lfloor{\frac{n}{2}}\right\rfloor.
\]

For further calculations, we will need explicit values of $\sin^{-1}\sqrt{x_j}$ and $\sinh^{-1}\sqrt{x_j}$, where the principal branches of the multivalued functions are implied. First, one can write
\[
\sin^{-1}\sqrt{x_j}=\xi_j-i\eta_j,\quad \sinh^{-1}\sqrt{x_j}=\varphi_j-i\psi_j,
\]
with $\xi_j,\eta_j,\varphi_j,\psi_j>0$.
Further, from elementary identities $1-2t=\cos(2\sin^{-1}\sqrt{t})$ and $1+2t=\cosh(2\sinh^{-1}\sqrt{t})$ one can see that 
\[
\cos(2\xi_j)\cosh(2\eta_j)=\cosh(2\varphi_j)\cos(2\psi_j)=1,
\]
\[
\sin(2\xi_j)\sinh(2\eta_j)=\sinh(2\varphi_j)\sin(2\psi_j)=\frac{\sin ^2\frac{\pi  (2 j-1)}{2 n}}{\cos\frac{\pi  (2 j-1)}{2 n}}.
\]
These equations can be easily solved to yield
\[
\xi_j=\psi_j=\frac{\pi (2 j-1)}{4 n},\quad \eta_j=\varphi_j=\frac12\sinh^{-1}\tan \frac{\pi  (2 j-1)}{2 n}.
\]
Thus
\begin{equation*}
    \sin^{-1}\sqrt{x_j}=\frac{\pi(2j-1)}{4n}-\frac{i}{2}\sinh ^{-1}\tan \frac{\pi  (2 j-1)}{2 n},
\end{equation*}
\begin{equation*}
    \sinh^{-1}\sqrt{x_j}=\frac{\pi(2j-1)}{4ni}+\frac{1}{2}\sinh ^{-1}\tan \frac{\pi  (2 j-1)}{2 n}.
\end{equation*}
Similarly
\[
\sin^{-1}\sqrt{y_j}=\frac{\pi(2j-1)}{4n}+\frac{i}{2}\sinh ^{-1}\tan \frac{\pi  (2 j-1)}{2 n},
\]
\[
\sinh^{-1}\sqrt{y_j}=-\frac{\pi(2j-1)}{4ni}+\frac{1}{2}\sinh ^{-1}\tan \frac{\pi  (2 j-1)}{2 n}.
\]
For $0\le s<\left\lfloor{\frac{n}{2}}\right\rfloor$ we have the partial fractions expansion
\[
\frac{t^{2s}}{Q_n(t)}=\sum_{j=1}^{n/2}\left(\frac{x_j^{2s}}{Q_n^{\,'}(x_j)}\frac{1}{t-x_j}+\frac{y_j^{2s}}{Q_n^{\,'}(y_j)}\frac{1}{t-y_j}\right).
\]
Calculations using the formulas above yield
\[
Q_n^{\,'}(x_j)=-Q_n^{\,'}(y_j)=\frac{4ni(-1)^j\cos^2\frac{\pi  (2 j-1)}{2 n}}{\sin\frac{\pi  (2 j-1)}{2 n}\left(1+\cos^2\frac{\pi  (2 j-1)}{2 n}\right)}.
\]
Now substitute this into the formula above.\qed

Using Lemma \ref{lemma5} and the following consequence of Lemma \ref{lemma2}
\[
\int\limits_0^1\frac{1}{4t^2+\frac{\sin ^4\frac{\pi  (2 j-1)}{2 n}}{\cos ^2\frac{\pi  (2 j-1)}{2 n}}}\frac{dt}{\sqrt{1-t^2}}=\frac{\pi}{2}\frac{\cot^2\frac{\pi  (2 j-1)}{2 n}}{1+\cos^2\frac{\pi  (2 j-1)}{2 n}},
\]
one can easily complete the proof of Theorem \ref{th3}.

\section{Proof of Theorem \ref{th4}}\label{proof4}

Both $\cos \left(2 u \sin ^{-1}\sqrt{t}\right)$ and $\cosh \left(2 u \sinh ^{-1}\sqrt{t}\right)$ are polynomials in $t$ when $u$ is an integer. Therefore, proceeding exactly as in the previous section one obtains that the integral in question is
\[
\frac{\pi }{n}\sum _{j=1}^{n/2}(-1)^{j-1}\tan\frac{\pi  (2 j-1)}{2 n}\frac{\cosh \left(u \sinh ^{-1}\tan \frac{\pi  (2 j-1)}{2 n}\right)}{\cosh \left(n \sinh ^{-1}\tan \frac{\pi  (2 j-1)}{2 n}\right)}\cos\frac{\pi  (2 j-1) u}{2 n}.
\]
Using 
\[
\frac{\cosh \left(uz\right)}{\cosh \left(n z\right)}=\frac{1}{n}\sum _{y=1}^n\frac{(-1)^{y-1}\sin\frac{\pi  (2 y-1)}{2 n}}{\cosh z-\cos\frac{\pi  (2 y-1)}{2 n}}\cos\frac{\pi  (2 y-1)u}{2 n}
\]
and the trivial identity $\cosh \left(\sinh ^{-1}\tan \frac{\pi  (2 j-1)}{2 n}\right)=1/\cos\frac{\pi  (2 j-1)}{2 n}$, one can rewrite this sum as a symmetric double sum
\[
\frac{\pi }{2n^2}\sum _{j,y=1}^{n}(-1)^{j+y}\frac{\sin\frac{\pi  (2 j-1)}{2 n}\sin\frac{\pi  (2 y-1)}{2 n}}{1-\cos \frac{\pi  (2 j-1)}{2 n}\cos \frac{\pi  (2 y-1)}{2 n}}\cos\frac{\pi  (2 j-1) u}{2 n}\cos\frac{\pi  (2 y-1) u}{2 n}.
\]
After some simple algebra the summand can be transformed as
\[
\frac{1}{2}(-1)^{j+y} \, \frac{\sin ^2\frac{\pi  (j+y-1)}{2 n}-\sin ^2\frac{\pi  (j-y)}{2 n}}{\sin ^2\frac{\pi  (j+y-1)}{2 n}+\sin ^2\frac{\pi  (j-y)}{2 n}}\left(\cos \frac{\pi (j+y-1)u}{n}+\cos\frac{\pi  (j-y)u}{n}   \right).
\]
$j-y$ and $j+y$ have the same parity. Taking into account this fact and periodicity of the summand, the summation can be performed independently over $j-y$ and $j+y$ to yield
\[
\frac{\pi }{4n^2}\sum _{x,s=1}^{n}\frac{\sin ^2\frac{\pi  (2x-1)}{2 n}-\sin ^2\frac{\pi  s}{n}}{\sin ^2\frac{\pi  (2x-1)}{2 n}+\sin ^2\frac{\pi  s}{n}}\left(\cos \frac{\pi (2x-1)u}{n}+\cos\frac{2\pi su}{n}   \right).
\]
Because of trivial identities
\[
\sum _{y=1}^{n} \cos\frac{\pi  (2 y+1)u}{n}=\sum _{x=1}^{n} \cos\frac{\pi  x u}{n}=0
\]
valid for integer $0<u<n$, and the summation formulas
\[
\sum _{x=1}^{n} \frac{\sinh ^2z}{\sinh ^2z+\sin ^2\frac{\pi  x}{n}}=\frac{n \coth (n z)}{\coth z},
\]
\[
\sum _{x=1}^{n} \frac{\cosh ^2z}{\sinh ^2z+\sin ^2\frac{\pi  (2x-1)}{2n}}=\frac{n \tanh (n z)}{\tanh z},
\]
the sum under consideration becomes
\[
\sum _{y=1}^n  \frac{\coth \left(n \sinh ^{-1}\sin \frac{\pi  (2 y-1)}{2 n}\right)}{\coth \left(\sinh ^{-1}\sin \frac{\pi  (2 y-1)}{2 n}\right)}\cos \frac{\pi (2 y-1)u}{n}+\sum _{x=1}^{n-1}  \frac{\sin ^2\frac{\pi  x}{n}}{1+\sin ^2\frac{\pi  x}{n}}\frac{\tanh \left(n \sinh ^{-1}\sin \frac{\pi  x}{n}\right)}{\tanh \left(\sinh ^{-1}\sin \frac{\pi  x}{n}\right)}\cos \frac{2\pi xu}{n}.
\]
Now it is easy to bring this to the form stated in the theorem.

\section{Discussion}\label{disc}

(i) A different proof of theorem \ref{th1} was given by P. Teruo Nagasava in a Math Stack Exchange post \cite{nagasava}. By using a clever substitution he was able to reduce the integral to an integral of a meromorphic function over the real line, and then use contour integration to evaluate it. From his solution one can also easily understand why the integral has the same value for all odd $n$.

(ii) When $u=0$, the sum in theorems \ref{th3} and \ref{th4} is simplified as
\[
\sum _{j=1}^{n/2} \frac{(-1)^{j-1} \tan \frac{\pi  (2 j-1)}{2 n}}{\cosh \left(n \sinh ^{-1}\tan \frac{\pi  (2 j-1)}{2 n}\right)}=\sum_{y=1}^n \frac{\coth \left(n \sinh ^{-1}\sin\frac{\pi  (2 y-1)}{2 n}\right)}{\coth \left(\sinh ^{-1}\sin\frac{\pi  (2 y-1)}{2 n}\right)}-\frac{n}{2}.
\]

(iii) Theorem \ref{th4} implies the integration formula
\begin{align*}
\frac{1}{\pi}\int\limits_{-1}^1&\frac{{\sin \left(2 n \sin ^{-1}\sqrt{t}\right)}}{\cos \left(2 n \sin ^{-1}\sqrt{t}\right)+\cosh \left(2 n \sinh ^{-1}\sqrt{t}\right)}\frac{\sqrt{t}\,dt}{\left(t-\sin ^2\frac{\pi  j}{2 n}\right)\sqrt{1+t}}\\&\phantom{........................}=1-\frac{ \sin\frac{\pi  j}{2 n} }{\sqrt{1+\sin ^2\frac{\pi j}{2 n}}}\left\{\tanh\left(n \sinh ^{-1}\sin\frac{\pi  j}{2 n}\right)\right\} ^{(-1)^j},
\end{align*}
where $j$ is an integer.

(iv) We briefly discuss the motivation behind the theorems presented in this paper. Let us introduce the notation
\begin{equation}\label{notation}
    \alpha_z=2n\sinh^{-1}\sin\frac{\pi z}{2n},
\end{equation}
where we assume the principal branches of the multivalued functions. With this definition one can rewrite the integral in Theorem \ref{th1} with $a=1$ as
\[
\int\limits_0^{1}\frac{\sin \bigl(n \sin^{-1}t\bigr)\sinh \bigl(n \sinh^{-1}t\bigr)}{\cos \bigl( 2 n \sin^{-1}t\bigr)+\cosh \bigl(2 n \sinh^{-1}t\bigr)}\frac{dt}{t \sqrt{1-t^4}}=\frac{\pi}{n}\int\limits_0^n\frac{\sin\frac{\pi x}{2}\sinh\frac{\alpha_x}{2}}{\cos\pi x+\cosh\alpha_x}\frac{dx}{\sinh\frac{\alpha_x}{n}}.
\]
As we will now show, the last integral has an interesting symmetry.

When $y$ is real, then $\alpha_{iy}$ is purely imaginary. Let us define $y_*$ by the equation
\[
\alpha_{iy_*}=\pi i n,
\]
and consider the integral over an interval on the imaginary axes
\[
J=\int\limits_{iy_*}^0\frac{\sin\frac{\pi z}{2}\sinh\frac{\alpha_z}{2}}{\cos\pi z+\cosh\alpha_z}\frac{dz}{\sinh\frac{\alpha_z}{n}}=\int\limits_0^{y_*}\frac{\sinh\frac{\pi y}{2}\sin\frac{\alpha_{iy}}{2i}}{\cos\pi y+\cos(\alpha_{iy}/i)}\frac{dy}{\sin\frac{\alpha_{iy}}{in}}.
\]
Making change of variables $\alpha_{iy}=\pi i s$ in \ref{notation} we get
\[
\sin\frac{\pi s}{2n}=\sinh\frac{\pi y}{2n},
\]
which implies that
\[
\pi y=\alpha_s.
\]
Since $\cos\frac{\pi s}{n}+\cosh\frac{\pi y}{n}=2$, it is easy to show that
\[
\frac{dy}{\sin\frac{\pi s}{n}}=\frac{ds}{\sinh\frac{\alpha_s}{n}}.
\]
Thus the integral under consideration becomes
\[
J=\int\limits_0^n\frac{\sin\frac{\pi s}{2}\sinh\frac{\alpha_s}{2}}{\cos\pi s+\cosh\alpha_s}\frac{ds}{\sinh\frac{\alpha_s}{n}}.
\]
To recap what we have just showed:
\[
\int\limits_{iy_*}^0\frac{\sin\frac{\pi z}{2}\sinh\frac{\alpha_z}{2}}{\cos\pi z+\cosh\alpha_z}\frac{dz}{\sinh\frac{\alpha_z}{n}}=\int\limits_0^n\frac{\sin\frac{\pi s}{2}\sinh\frac{\alpha_s}{2}}{\cos\pi s+\cosh\alpha_s}\frac{ds}{\sinh\frac{\alpha_s}{n}},\qquad y_*=\frac{2n}{\pi}\ln(1+\sqrt{2}).
\]
In words, the integral of the function
\[
\frac{\sin\frac{\pi z}{2}\sinh\frac{\alpha_z}{2}}{\cos\pi z+\cosh\alpha_z}
\]
taken over a certain segment of the imaginary axis, turns out to be equal to the integral of this function taken over a segment of the real axis.

The integral in theorem \ref{th1} has been chosen to have the same kind of symmetry:
\[
\int\limits_0^{1}\frac{t^{2s}}{\cos (2 n \sin^{-1}\!\sqrt{t})+\cosh \left(2 n \sinh^{-1}\!\sqrt{t}\right)}\frac{dt}{\sqrt{1-t^2} }=\frac{\pi}{n}\int\limits_0^n\frac{\left(\sin\frac{\pi x}{2n}\right)^{4s+2}}{\cos\pi x+\cosh\alpha_x}\frac{dx}{\sinh\frac{\alpha_x}{n}},
\]
\[
\int\limits_{iy_*}^0 \frac{\left(\sin\frac{\pi z}{2n}\right)^{4s+2}}{\cos\pi z+\cosh\alpha_z}\frac{dz}{\sinh\frac{\alpha_z}{n}}=\int\limits_0^n\frac{\left(\sin\frac{\pi x}{2n}\right)^{4s+2}}{\cos\pi x+\cosh\alpha_x}\frac{dx}{\sinh\frac{\alpha_x}{n}},\qquad y_*=\frac{2n}{\pi}\ln(1+\sqrt{2}).
\]

\section{Some open questions}\label{open}

First question is concerned with an extension of theorem \ref{th4} to the non-symmetric case $a\neq 1$.

The second question is concerned with some integration formulas similar to \ref{glaisher1} and \ref{glaisher2}. By contour integration it is fairly easy to prove the integration formula
\begin{equation}\label{glaisher3}
    \int\limits_0^\infty{\mathrm{Im}}\left\{{\frac{1}{\cosh \left(e^{-i \alpha} x\right) \cosh \left(e^{-i \beta} x\right)}}\right\}\frac{dx}{x}=\frac{\alpha+\beta}{2}.
\end{equation}
For this, we take the function $\left\{{\cosh \left(e^{i \alpha} z\right) \cosh \left(e^{i \beta} z\right)}\right\}^{-1}$ and integrate it along the contour composed of two rays $\arg z=0$, $\arg z=\alpha+\beta$, and two circular arcs, one around the origin, and the other around the complex infinity. The integrals over the rays combine to the integral in \ref{glaisher3} multipled by $2i$. The integral over the arc at infinity vanishes, while the arc around the origin gives the contribution $-i(\alpha+\beta)$. The final step is to notice that there are no poles inside the contour.

Writing out the integrand explicitly and taking the sum of two integrals with $\pm\beta$ one obtains
\[
\int\limits_0^\infty\frac{\sin(x\sin \alpha)\sinh(x \cos \alpha)\cos(x\sin \beta)\cosh(x \cos \beta)}{\left\{\cosh(2x\cos \alpha)+\cos (2x\sin \alpha)\right\}\left\{\cosh(2x\cos \beta)+\cos (2x\sin \beta)\right\}}\frac{dx}{x}=\frac{\alpha}{8}.
\]
Although this formula was derived for real $\alpha$ and $\beta$ it can be continued analytically to complex values. The question is to find a finite analog of this integral similar to the one in theorem \ref{th1}. Direct naive extension does not work, even for the integral with $\alpha=\beta$
\[
\int\limits_0^\infty\frac{\sin(x\sin \alpha)\sinh(x \cos \alpha)}{\left\{\cosh(x\cos \alpha)+\cos (x\sin \alpha)\right\}^2}\frac{dx}{x}=\frac{\alpha}{2}.
\]
It might be worth mentioning the sum closely related to these integrals
\[
\sum_{n=1}^\infty \frac{\chi(n)}{n}\cdot\frac{ \cos \left(\frac{\pi n\cos \theta}{2}\right) \cosh \left(\frac{\pi n\sin \theta}{2}\right)}{\cos (\pi n\cos \theta)+\cosh (\pi n \sin \theta)}=\frac{\pi}{16},
\]
where $\chi(n)=\sin\frac{\pi n}{2}$ is Dirichlet character modulo $4$.

\end{document}